\documentclass{article}
\usepackage{graphicx}
\usepackage[cmex10]{amsmath}
\usepackage{amsthm}

\usepackage{amssymb, subfigure,float, url, color}
\usepackage[colorlinks=true]{hyperref}
\usepackage{pdfsync}

\newcommand{\condition}{$(\Omega)$}

\newcommand{\g} {\gamma}

\newcommand{\hh} {{\cal H}}
\newcommand{\spann} {\operatorname{span }}

\newcommand{\N} {\mathbb{N}}

\newcommand{\R} {\mathbb{R}}
\newcommand{\con}{{\cal C}}

\renewcommand{\geq}{\geqslant}
\renewcommand{\leq}{\leqslant}
\newcommand{\Lie}{\operatorname{Lie}}
\newcommand{\ocpe}{(\mathrm{OCP})_{\varepsilon}}
\newcommand{\ocp}{(\mathrm{OCP})}
\newcommand{\ocpse}{(\mathrm{OCPS})_{\varepsilon}}
\newcommand{\ocps}{(\mathrm{OCPS})}

\newtheorem{theorem}{Theorem}  
\newtheorem{proposition}[theorem]{Proposition}

\newtheorem{definition}{Definition}
\newtheorem{lemma}[theorem]{Lemma}
\newtheorem{example}{Example}
\theoremstyle{definition}\newtheorem{remark}{Remark}

\newcommand{\bqnn}{\begin{equation*}}
\newcommand{\eqnn}{\end{equation*}}

 \title{\LARGE \bf
Regularization of chattering phenomena via bounded variation controls}
\author{
Marco~Caponigro\footnote{Conservatoire National des Arts et M\'etiers, Equipe M2N, 75003 Paris, France ({\tt marco.caponigro@cnam.fr}).} 
\quad 
Roberta~Ghezzi\footnote{Institut de Math\'ematiques de Bourgogne, Universit\'e de Bourgogne-Franche Comt\'e, Dijon, France ({\tt roberta.ghezzi@u-bourgogne.fr}).}
\quad
Benedetto Piccoli\footnote{
Department of Mathematical Sciences and Center for Computational and Integrative Biology, Rutgers University, Camden, NJ 08102, USA ({\tt piccoli@camden.rutgers.edu}).}
\quad
Emmanuel Tr\'elat\footnote{
Sorbonne Universit\'es, UPMC Univ Paris 06, CNRS UMR 7598, Laboratoire Jacques-Louis Lions, Institut Universitaire de France, F-75005, Paris, France ({\tt emmanuel.trelat@upmc.fr).}
 }}

\date{}

\begin{document}


\maketitle

\begin{abstract}
In control theory, the term \emph{chattering} is used to  refer to  strong oscillations of controls, such as an infinite number of switchings over a compact interval of times. In this paper we focus on three typical occurences of chattering: the Fuller phenomenon, referring to situations where an optimal control switches an infinite number of times over a compact set; the Robbins phenomenon, concerning optimal control problems with state constraints, meaning that the optimal trajectory touches the boundary of the constraint set an infinite number of times over a compact time interval; the Zeno phenomenon, referring as well to an infinite number of switchings over a compact set, for hybrid optimal control problems.
From the practical point of view, when trying to compute an optimal trajectory,  for instance by means of a shooting method,  chattering  may be a serious obstacle to convergence.

In this paper we propose a general regularization procedure, by adding an appropriate penalization of the total variation.
This produces a quasi-optimal control, and we prove that the family of quasi-optimal solutions converges to the optimal solution of the initial problem as the penalization tends to zero. Under additional assumptions, we also quantify the quasi-optimality property by determining a speed of convergence of the costs.
\end{abstract}



\section{Introduction}
Chattering phenomena in optimal control have been known since the first example presented in~\cite{fuller}.
Roughly speaking, chattering refers to strong oscillations of the optimal control switching infinitely many times over a finite time interval. To explain this behavior, let us recall the famous example in \cite{fuller}, also known as {\it Fuller's phenomenon}. Given $T>0$ arbitrary, consider the control system in $\R^2$
\begin{equation}\label{intro1}
\dot x_1=x_2,\quad \dot x_2=u,
\end{equation}
with controls $u:[0,T]\to [-1,1]$, and consider the optimal control problem consisting of minimizing the cost functional
\begin{equation}\label{intro2}
\int_0^T x_1^2(t)\,dt ,
\end{equation}
over all trajectories of \eqref{intro1} steering an (arbitrary) initial point $(x_1^0,x_2^0)$ to the origin, i.e., such that
$$
x_1(0)=x_1^0,\quad x_2 (0)=  x_2^0,
\qquad
x_1(T)= 0,\quad x_2(T)=0.
$$
As is well known, there exists a unique optimal control $u:[0,T]\to[-1,1]$, satisfying
$$
u(t)= 
\begin{cases}
1, & t\in(t_{2k},t_{2k+1}), \quad  k\in\N,\\
-1, & t\in(t_{2k+1},t_{2k+2}), \quad  k\in\N,
\end{cases}
$$
where $(t_k)_{k\in\N}$ is an increasing sequence of switching times, depending on the initial condition $(x_1^0,  x_2^0)$ and converging to $T$.
Although, at the first sight, one could think that this strong oscillation property is a kind of aberration due to specific symmetries of the system, it turns out that this chattering behavior is rather typical. Indeed, it was later shown in \cite{kupka} that the set of single-input optimal control problems which have a control-affine Hamiltonian and whose solution is chattering is an open semi-algebraic set (see also~\cite{borisov}), showing therefore that chattering is a common phenomenon in optimal control.

Controls enjoying a chattering property have been found for a variety of problems: besides the ones mentioned above,  a similar phenomenon also concerns state-constrained problems and hybrid systems. In~\cite{robbins}, Robbins studied an optimal control problem with an inequality state constraint of third order, and he showed that the optimal trajectory touches the constraint's boundary at an infinite sequence of isolated points converging to a point at the boundary (and however the optimal control has finite total variation).
In the framework of hybrid systems, chattering is often called Zeno phenomenon and is due to  trajectories whose discrete part jumps infinitely many times over a finite time interval (see, for instance, the examples in~\cite{ZJLS01}).

Although chattering cannot be considered as a degeneracy phenomenon (see \cite{Schattler_book}), chattering may however cause some difficulties in theoretical and numerical aspects of optimal control. 

From the theoretical point of view, due to the lack of a positive length interval where the control function is continuous when chattering occurs, finding  necessary and sufficient optimality conditions becomes much more intricate (see \cite{bonnans-hermant} for state-constrained problems). Some results in this sense were proved in~\cite{borisov}, yet the problem is not completely understood in other contexts, such as state-constrained problems or hybrid systems  \cite{zeno-numpb}. Another delicate issue comes from the  study of   regularity  properties of optimal syntheses \cite{piccbosc,picsuss}.

From the numerical point of view, chattering phenomena may be an obstacle to the convergence of numerical methods applied to optimal control problems, in particular when using indirect methods.
Indeed, chattering implies ill-posedness of shooting methods (non-invertible Jacobian) \cite{bonnans-interior,bonnans-wellposedness}.
When chattering occurs in an optimal control problem, it is therefore required to develop an adequate numerical method in order to compute a good approximation of the optimal control. This problem has been raised in \cite{ZTC2,ZTC} for the optimal control of the attitude of a launcher, in which chattering may occur, depending on the terminal conditions under consideration. After having observed that chattering was indeed causing the failure of the shooting method, the authors have proposed two remedies: one is based on a specific homotopy combined with the shooting method, and the other consists of using a direct method with a finite number of arcs. However, on the one part these remedies remain specific to the problem studied thereof, and on the other part there is no convergence result that would show and quantify the quasi-optimality property.

\medskip

In this paper, we propose a general regularization procedure, consisting of penalizing the cost functional with a total variation term. Our approach is valid for general classes of nonlinear optimal control problems. For a bang-bang scalar control, the total variation of the control is {proportional} to the number of switchings. In the case where the Fuller phenomenon occurs, the total variation is infinite.
Hence, with such a penalization term, the optimal control does not chatter, and its numerical computation is then a priori feasible. 
The total variation term is a penalization that regularizes the optimal control. Under appropriate assumptions of small-time local controllability, we prove in Theorem~\ref{thmtv} that, if the weight $\varepsilon$ of the total variation term in the cost functional tends to zero, then the regularized optimal control problem $\Gamma$-converges to the initial optimal control problem, meaning that the optimal cost and the optimal solution of the regularized problem converge respectively to the optimal cost and the optimal solution of the initial problem, as the parameter $\varepsilon$ tends to zero. 
This shows that, when this total variation regularization is used, the optimal control that one may then compute numerically is quasi-optimal, with a good rate of optimality.

In order to quantify quasi-optimality, it remains then to determine at what speed the cost of the regularized problem converges to the cost of the initial problem, as the weight of the total variation term tends to zero. This can be done by estimating explicitly the rate of convergence of the cost along suboptimal regimes obtained by suitable truncations of the chattering one in terms of switching times. In the existing literature, such results, related to truncation, were obtained in~\cite{zelman08,zelzel} for  small perturbations of the Fuller's problem. In those papers, the authors exhibited a sequence of suboptimal regimes for the specific optimal control problem \eqref{intro1}-\eqref{intro2}, and they proved that the cost converges with the same rate as the sequence of switching times (of the chattering control). 
Our Theorem~\ref{thm:rate} establishes a polynomial rate of convergence for the costs, for general nonlinear optimal control problems, under appropriate controllability assumptions, and under the additional assumption that the time-optimal map is H\"older continuous.    
Note that, for the specific case considered in~\cite{zelzel}, the rate of convergence is exponential as a function of the number of switchings. Likewise, for the class of systems considered in \cite{kupka}, the switching times converge exponentially to the final time. Whether a slower rate of  convergence is ``typical'' remains an open question. 

Finally, we treat by total variation regularization two other general cases where chattering occurs:
\begin{itemize}
\item For optimal control problems involving state-constraints, under adequate controllability assumptions, Theorem \ref{robbinseps} provides a regularization result for optimal trajectories having an infinite sequence of contact points with the constraint's boundary (Robbins phenomenon). Here, the penalization term essentially counts the contact points with the constraint's boundary.
\item For hybrid optimal control problems, Theorem \ref{thmtvzeno} provides a convergence result to regularize the Zeno phenomenon, obtaining estimates of the cost convergence as the number of location switchings grows.
\end{itemize}

The paper is organized as follows. In Section~\ref{sec:main} we present the main results on the regularization by total variation penalization of chattering phenomena (Fuller, Robbins and Zeno). Section~\ref{sec:proof} is devoted to prove the main results.
In Appendix~\ref{sec:remarks} we provide some additional results concerning the controllability condition required in Theorem~\ref{thm:rate}. Finally, we provide in Appendix~\ref{sec:app} an existence result for optimal control problems having a total variation term in the cost functional, without any convexity assumptions.

\section{Main results}\label{sec:main}
\subsection{Regularization of the Fuller phenomenon}\label{subsec_Fuller}
Let $N$ and $m$ be nonzero integers. Consider the control system
\begin{equation}\label{dinamica}
\dot x=f(x,u),\quad u\in{\cal U}, \tag{$\Sigma$}
\end{equation}
where 
$f\in\con^\infty(\R^N\times\R^m,\R^N)$, $f(0,0)=0$, and
\begin{equation}\label{defcalU}
{\cal U}=\{u(\cdot) \mbox{ measurable} \mid u(t) \in {\bf U} \mbox{ for a.e. } t \},
\end{equation}
with   ${\bf U}\subset \R^m$ a measurable subset {containing $0$}.
Denote by 
$$
{\cal F}=\{f(\cdot, u):\R^N\to\R^N\mid  u\in{\bf U}\},
$$ 
the family of vector fields associated with the dynamics of~\eqref{dinamica}.

A control $u \in \mathcal{U}$ is called \emph{admissible} if it steers the system~\eqref{dinamica} from a given (arbitrary) initial point to the origin in finite time denoted $t(u)$.

Given an initial state $x_0\in\R^N$, a function $L\in\con^0(\R\times\R^N\times\R^m)$ (called Lagrangian), we consider the optimal control problem 
\begin{equation*}
\begin{cases}
&\displaystyle{\min_{u\in{\cal U}}  \int_0^{t(u)}L(s,x(s),u(s))\,ds ,} \\
&\dot x=f(x,u),\quad u\in{\cal U},\\
&x(0)=x_0,\quad x(t(u))=0.
\end{cases}
\eqno{(\mathrm{OCP})}
\end{equation*}
The final time in $\ocp$ may be fixed or free. If it is fixed to some $T>0$, then of course one has to replace $t(u)$ with $T$ everywhere.

Throughout the section, we make the following assumptions:
\begin{itemize}
\item for every $(t,x)\in\R\times\R^N$, the set
\begin{equation}\label{defV}
V(t,x)=\{(f(x,u),L(t,x,u)+\gamma)\mid u \in{\bf U}, \gamma \geq 0\}
\end{equation}
is convex;
\item ${\bf U}$ is compact, and there exists $b>0$ such that, for every admissible control $u\in{\cal U}$, we have
\begin{equation}\label{hypub}
t(u)+\Vert x_u(\cdot)\Vert_{\infty}\leq b .
\end{equation}
\end{itemize}
The first assumption means that the epigraph of extended velocities is convex. It is satisfied, for example, for control-affine systems with control-affine or quadratic cost.

These are classical assumptions used to derive existence results (see, for instance, \cite{cesari-book, benlibro, trelat_book}). 
Under these assumptions, the optimal control problem $\ocp$ has at least one optimal  solution $x^*(\cdot)$, associated with a control $u^*:[0,t(u^*)]\to{\bf U}$.

\medskip

It may occur that the control $u^*$ chatters. In this case, as discussed in the introduction, this may cause the failure of numerical methods to compute it.
To overcome this problem, we next propose a regularization of the optimal control problem $\ocp$ by adding to the cost functional a total variation term, penalizing oscillations, with a small weight $\varepsilon$.

Given any $\varepsilon\geq 0$, we consider the optimal control problem 
\begin{equation*}\label{costoeps}
\begin{cases}
&\displaystyle{\min_{u\in{\cal U}}  \left(\int_0^{t(u)}L(s,x(s),u(s))\,ds + \varepsilon\, \mathrm{TV}(u)\right),} \\
&\dot x=f(x,u),\quad u\in{\cal U},\\
&x(0)=x_0,\quad x(t(u))=0.
\end{cases}
\eqno{(\mathrm{OCP})_{\varepsilon}}
\end{equation*}
Here, $\mathrm{TV}(u)$ designates the total variation of the function $u:[0,t(u)] \to \R^{m}$, and it is defined by
$$
\mathrm{TV}(u) = \sup \sum_{i=1}^{p} \|u(t_{i}) - u(t_{i-1})\|,
$$ 
the supremum being taken over all possible partitions $0=t_{0} < t_{1}< \cdots < t_{p} =t(u)$ of the interval $[0,t(u)]$. For instance, if $m=1$ and if $u$ is a  piecewise constant function taking values in $\{0,1\}$, then $\mathrm{TV}(u)$ is simply equal to the number of switchings.  
A function $u:[0,t(u)] \to \R^{m}$ is said to have bounded variation if $\mathrm{TV}(u)<+\infty$.

The rationale for introducing the term $\varepsilon\, \mathrm{TV}(u)$ in the cost of~\eqref{costoeps} is to penalize highly oscillating controls in order to avoid \emph{chattering} in the sense of Definition~\ref{def:chat}
below.

\begin{definition}\label{def:chat}
By \emph{chattering} control we mean a measurable function $u:[0,t(u)]\to{\bf U}$ such that
there exists an increasing sequence $\{t_n\}_{n\in \N}$ converging to $t(u)$ with the property that
$\mathrm{TV}(u|_{[0,t_{n}]}) < +\infty$
for every $n\in\N$,
and 
$$
\lim_{n\to+\infty}\mathrm{TV}(u|_{[0,t_n]})=+\infty.
$$ 
\end{definition}

The optimal control problem $\ocpe$ is seen as a regularization of $\ocp$. We are next going to prove that any optimal solution $\ocpe$ converges {uniformly} to an optimal solution of $\ocp$, thus providing a quasi-optimal solution that does not chatter.

Recall that the control system~\eqref{dinamica} is \emph{small-time locally controllable} (STLC) at $x_{0} \in \R^{N}$ if, for every $\delta >0$, there exists a neighborhood ${\cal N}_{\delta}$ of $x_{0}$ such that every $x_{1} \in {\cal N}_{\delta}$ can be reached by $x_{0}$ within time $\delta$ with a control $u \in \mathcal{U}$.

In the sequel,
$\Lie_x{\cal F}$ denotes the vector field Lie algebra generated by ${\cal F}$ evaluated at $x$, that is $\Lie_x{\cal F} = \{V(x)\ \mid\ V \in \Lie{\cal F}\}$, where $\Lie{\cal F}  = \mathrm{span}
\{[f_{1},[ \dots[f_{k+1},f_{k}]\ldots]] \mid  f_{i} \in {\cal F}, k \in \N\}$.
 
\begin{theorem}\label{thmtv}
Assume that $\Lie_0{\cal F}=\R^N$, and the control system~\eqref{dinamica} is small-time locally controllable at $0$.
Then, for every $\varepsilon >0$, the optimal control problem $\ocpe$ 
has at least one solution. Moreover, for every optimal solution $x_\varepsilon(\cdot)$ of $\ocpe$, associated with a control $u_\varepsilon:[0,t(u_\varepsilon)]\to {\bf U}$, we have
\begin{equation}\label{convcost}
\lim_{\varepsilon\to 0} \int_0^{t(u_\varepsilon)} L(t,x_\varepsilon(t),u_\varepsilon(t))\,dt=\int_0^{t(u^*)}L(t,x^*(t),u^*(t))\,dt ,
\end{equation}
and $x_\varepsilon(\cdot)$ converges uniformly to an optimal solution of $\ocp$.
\end{theorem}

Theorem~\ref{thmtv} establishes the existence of a non-chattering control
$u_\varepsilon$ which is quasi-optimal for $\ocp$ in the sense that the cost of $u_\varepsilon$ converges to the optimal value of $\ocp$.

\begin{remark}\label{rk:1706}
The Lie algebra and small-time controllability assumptions, although generic, may be slightly weakened without altering the conclusion of the theorem: they can be replaced by assuming local controllability in a neighborhood of the origin in arbitrarily small time and with piecewise constant controls. The fact that the latter assumption is weaker follows from a well known result due to Krener (see for instance \cite[Corollary~8.3]{agrachev-book}). 
\end{remark}

\begin{remark}
Note that we have assumed $f$ to be smooth, in order to give a sense to Lie brackets (we assume that $\Lie_0{\cal F}=\R^N$ in the theorem). In contrast, we only need the Lagrangian function $L$ to be continuous. Besides, $L$ may depend on $t$ but it is important that the dynamics $f$ is autonomous (the fact that $f(0,0)=0$ is useful in the proofs).
\end{remark}

In the case where the optimal control of $\ocp$ chatters and therefore cannot be computed by means of a shooting method, the total variation term in $\ocpe$ plays the role of a regularization, and the control $u_\varepsilon$ does not chatter and can be computed numerically.
Theorem~\ref{thmtv} establishes that $u_\varepsilon$ is quasi-optimal, and hence it is reasonable to replace $\ocp$ by $\ocpe$ when chattering occurs, in order to ensure the convergence of a shooting method.

Theorem~\ref{thmtv} establishes the convergence \eqref{convcost} of the costs. It is then interesting to derive a speed of convergence. This is possible under additional assumptions, as we are going to see next.

We need the following ``strong'' notion of controllability which requires a uniform bound on the total variation of the control and a steering time comparable with the minimum time. To this purpose, we define the time-optimal map $x_{0}\mapsto \Upsilon(x_{0})$ associated with the control system~\eqref{dinamica}, by
\begin{equation}\label{iu}
\Upsilon(x_{0}) = \inf \{t >0 \mid \dot x = f(x,u),\ x(0)=x_{0},\ x(t)=0 \}.
\end{equation}  

\begin{definition}\label{def:omega}
We say that the control system \eqref{dinamica} satisfies \condition\ at $0$ if
\begin{itemize}
\item[$(\Omega_{1})$] the control system \eqref{dinamica} is STLC at $0$;
\item[$(\Omega_{2})$] there exist a neighborhood $\cal N$ of $0$ and $M>0$ such that, for every $y \in \mathcal{N}$, there exists $u:[0,\tau_{y}] \to {\bf U}$ such that
  $u$ steers $y$ to $0$ in time  $\tau_{y}$,
  $\tau_{y} \leq M \Upsilon(y)$,
$\mathrm{TV}(u) \leq M$.
\end{itemize}
\end{definition}

We provide in Appendix~\ref{sec:remarks} some comments on Definition~\ref{def:omega} and some results on the relationships between the properties \condition, STLC, and the regularity of $\Upsilon$.

In the sequel, ${\cal C}^{0,\alpha}$ designates the class of H\"older continuous functions with exponent $\alpha$.

\begin{theorem}\label{thm:rate}
Assume that:
\begin{itemize}
\item[$(i)$] the control system \eqref{dinamica} satisfies \condition\ at $0$;
\item[$(ii)$] the optimal control $u^*$ of $\ocp$ either has bounded total variation, or is chattering and its sequence of switching times $(t_{n})_{n\in\N}$ satisfies $(t(u^*)-t_{n}) = \mathrm{O}(n^{-\beta})$ for some $\beta>0$;
 \item[$(iii)$] the time-optimal map is $\con^{0,\alpha}$ for some $\alpha \in (0,1]$ in a neighborhood   
of $0$.  
\end{itemize}
Then, for every $\varepsilon>0$, the optimal control problem $\ocpe$ has at least one solution. Moreover, for every optimal solution $x_{\varepsilon}(\cdot)$ of $\ocpe$, associated with a control $u_{\varepsilon}:[0,t(u_\varepsilon)]\to {\bf U}$, we have 
\begin{align}
 \int_0^{t(u_\varepsilon)}L(t,x_{\varepsilon}(t),u_{\varepsilon}(t))\,dt-  \int_0^{t(u^*)} & L(t,x^*(t),u^*(t))\,dt  = \nonumber\\  &  \quad =\mathrm{O}\left( \varepsilon^{\frac{\alpha\beta}{1+\alpha\beta}} \right). \label{convcost2}
\end{align}
\end{theorem}

\begin{remark}
For linear control systems and for driftless control-affine systems,  \condition\ is related to controllability. Sufficient conditions guaranteeing that \condition\ holds true can be found in~\cite{sharon,suss-bound} for single-input control systems. For more general control-affine systems, \condition\ is related to the \emph{Exact State Space Linearizability Problem} (see Appendix~\ref{sec:remarks}).
\end{remark}

\begin{remark}
Assumption $(ii)$ is verified for a large class of systems having an exponential rate of accumulation of switchings (see \cite{kupka}). In this case the convergence rate is $\mathrm{O}(\varepsilon^{\gamma})$ for every $\gamma<1$.
\end{remark}

\begin{remark}
Sufficient conditions for Assumption $(iii)$ have been established in~\cite[Theorem~3.3, 3.10, 3.12]{stefani-holder}, where  the authors provide an estimate on the H\"older exponent. 
\end{remark}

\subsection{Regularization of the Robbins phenomenon for problems with state constraints}\label{ro}
In this section, we consider the general optimal control problem $\ocp$ of the previous section with additional  state constraints.
Namely, let 
\begin{equation}\label{state}
\con=\{x\in\R^N \mid h_1(x)\geq 0,\dots, h_l(x)\geq 0\},
\end{equation}
where $h_1, \ldots, h_{l}$ are continuous functions and consider the optimal control problem
\begin{equation*} 
\begin{cases}
&\displaystyle{\min_{u\in{\cal U}} \int_0^{t(u)}L(s,x(s),u(s))\,ds ,} \\
&\dot x=f(x,u),\quad u\in{\cal U},\\
& x(t)\in\con,\quad t\in[0,t(u)], \\ 
&x(0)=x_{0},\quad x(t(u))=0.
\end{cases}
\eqno{(\mathrm{OCPS})}
\end{equation*}
The notations and the assumptions on the dynamics are the same as in Section \ref{subsec_Fuller}. In particular, we assume that the epigraph of extended velocities \eqref{defV} is convex, that ${\bf U}$ is compact, that \eqref{hypub} holds true, and that there exists at least one admissible trajectory satisfying the constraints.
Under these assumptions, the optimal control problem $\ocps$ has at least one optimal  solution $x^*(\cdot)$, associated with a control $u^*:[0,t(u^*)]\to{\bf U}$ (see \cite{cesari-book, benlibro, trelat_book}).

In \cite{robbins}, an instance of $\ocps$ is provided where the final point $0$ lies on the boundary $\partial\con$,    the solution $u^{*}$ is $\con^{1}$-smooth and  the trajectory $x^{*}(\cdot)$ corresponding to $u^{*}$ touches $\partial\con$ at a sequence of isolated points converging to the final point $0$.  In other words, the optimal trajectory is a concatenation of an infinite number of arcs contained in the interior of $\con$ and accumulating at the final point. We call this phenomenon the Robbins phenomenon.

To regularize this chattering effect, one needs to find suboptimal controls whose trajectories touch $\partial\con$ on a finite set. 
Introducing the total variation of the control as a penalization term, as it was done previously, does however not suffice to prevent the solution of the regularized problem from possibly intersecting $\partial\con$ infinitely many times. We next design a penalization term that rather counts the number of contact points with $\partial\con$.

Let ${\bf 1}_{\partial\con}:\R^N\to\{0,1\}$ be the indicator function of $\partial\con$, defined by
$$
{\bf 1}_{\partial\con}(x)=
\begin{cases}
1& \textrm{if}\ x\in\partial\con,\\
0& \textrm{if}\ x\notin\partial\con.
\end{cases}
$$  
Given an admissible control $u:[0,t(u)]\to {\bf U}$ with corresponding trajectory $x(\cdot)$, we define the function $X_u:[0,t(u)]\to\{0,1\}$ by 
$$
X_u(t) = {\bf 1}_{\partial\con}(x(t)) .
$$
For every $\varepsilon\geq 0$, we consider the optimal control problem 
\begin{equation*}\label{costoepsvinc}
\begin{cases}
&\displaystyle{\min_{u\in{\cal U}} \left( \int_0^{t(u)}L(s,x(s),u(s))\,ds +\varepsilon \, \mathrm{TV}(X_u) \right),} \\
&\dot x=f(x,u),\quad u\in{\cal U},\\
& x(t)\in\con, \quad t\in[0,t(u)],\\ 
&x(0)=x_{0},\quad x(t(u))=0.
\end{cases}
\eqno{(\mathrm{OCPS})_{\varepsilon}}
\end{equation*}

In the sequel, we consider the reachable set from $0$ with trajectories lying in the interior $\mathring{\mathcal{C}}$ of the constraint set $\mathcal{C}$ defined by \eqref{state}: let $\mathcal{A}^\mathcal{C}(0,(0,\delta),f)$ be the set of points accessible  from $0$ in time $t\in (0,\delta)$ by trajectories $x(\cdot)$ of the control system \eqref{dinamica} such that $x(t) \in \mathring{\mathcal{C}}$ for every $t \in (0,\delta)$.

\begin{theorem}\label{robbinseps}
Assume that $0\in\partial\con$ and that:
\begin{itemize}
\item[$(i)$] for every $\delta>0$, there exists a neighborhood $\mathcal{N}$ of $0$ such that
$\mathcal{N}\cap \mathring{\con} \subset \mathcal{A}^\mathcal{C}(0,(0,\delta),-f)$;
\item[$(ii)$] there exists a sequence of times $t_n$ converging to $t(u^*)$, with 
$$x^*([0,t(u^*)])\cap\partial\con\ \subset\ \{0\}\cup\{x^*(t_n)\ \mid\ n\in \N\};$$
\end{itemize}
Then, for every $\varepsilon >0$, the optimal control problem $\ocpse$ has at least one solution. Moreover, for every optimal solution $x_{\varepsilon}(\cdot)$ of $\ocpse$, associated with a control $u_\varepsilon:[0,t(u_\varepsilon)]\to {\bf U}$, we have
\begin{equation}\label{convcostvinc}
\lim_{\varepsilon\to 0} \int_0^{t(u_\varepsilon)}L(t,x_\varepsilon(t),u_\varepsilon(t))\,dt=\int_0^{t(u^*)}L(t,x^*(t),u^*(t))\,dt,
\end{equation}
and $x_\varepsilon(\cdot)$ converges uniformly to an optimal solution of $\ocps$.
\end{theorem}

\begin{remark}
Assumption $(i)$ is an adaptation of the classical small-time local attainability (STLA) property (see~\cite{krasta2001,marigonda}), but we require here that the admissible trajectories stay in the interior of the constraint $\mathcal{C}$. Hence Assumption $(i)$ may be seen as a generalization to nonlinear systems of the notion of small-time controllability with respect to a cone. 
Controllability with respect to a cone has been studied for linear control systems in \cite{krasta2008} (see also \cite{broucke2006}).
\end{remark}

\subsection{Regularization of the Zeno phenomenon for hybrid problems}\label{sec:zen}
In this section, we adapt Theorem~\ref{thmtv} to hybrid optimal control problems, where the dynamics involves a continuous and a discrete part. Let us first recall some basic notions on hybrid systems, without control (see, e.g., \cite{zeno-numpb}).
A \emph{hybrid system} is a collection $\hh=(Q,X,f,E,G,R)$ where
\begin{itemize}
\item  $Q$ is a finite set;
\item $X=\{X_q\}_{q\in Q}$ is a collection of subsets $X_q\subset\R^N$ called \emph{locations};
\item $f=\{f_q\}_{q\in Q}$ is a collection of smooth vector fields $f_q:\R^N\to \R^N$;
\item $E\subset Q\times Q$ is a subset of edges;
\item $G$ maps an edge $(q,q')\in E$ to a subset $G(q,q')\subset X_q$ called \emph{guard set};
\item $R$ maps a pair $((q,q'), x)\in E\times X_q$ to a subset $R((q,q'),x)\subset X_{q'}$. 
\end{itemize}
A \emph{trajectory} (or \emph{execution}) of $\hh$ is a triple $(\tau, q(\cdot), x(\cdot))$, where
\begin{itemize}
\item $\tau=\{\tau_i\}_{i=0}^M$ is a sequence of increasing positive numbers such that $\tau_0=0$ and $M\leq \infty$. We set $I=[0,\tau_M]$ if $M<+\infty$, $I=[0,\tau_M)$ if $M=\infty$;
\item $q:I\to Q$ is such that $q(t) = q_i$ constant on $[\tau_i,\tau_{i+1})$ for every $i=0,\dots M-1$;
\item for every $i=0,\dots, M-1$,  $x_i(\cdot)=x|_{(\tau_i,\tau_{i+1})}$ is an absolutely continuous function in $(\tau_i,\tau_{i+1})$, which can be continuously extended to $[\tau_i,\tau_{i+1}]$, and such that $x_i(t)\in X_{q_i}$;
\item for almost every $t\in(\tau_i,\tau_{i+1})$,
\begin{equation}\label{csh}
\dot x_i=f_{q_{i}}(x_i);
\end{equation}
\item for every $i=0,\dots, M-1$, one has $(q_i,q_{i+1})\in E$ and $x_i(\tau_{i+1})\in G(q_i,q_{i+1})$ and, for every $i=0,\dots,M-2$, one has $x_{i+1}(\tau_{i+1})\in R((q_i,q_{i+1}),x_i(\tau_{i+1}))$.
\end{itemize}
We say that $(\tau,q(\cdot), x(\cdot))$ is a \emph{Zeno trajectory} if $M=+\infty$ and $\tau_\infty<+\infty$. 

Given a hybrid system $\hh$, a \emph{Lagrangian} for $\hh$ is a family 
$L=\{L_q\}_{q\in Q}$, with $L_q:\R\times X_q\to\R$ such that, for every trajectory $(t,q(\cdot),x(\cdot))$ of $\hh$ and every $i=0,\dots, M-1$, the function $t\mapsto L_{q_i}(t, x_i(t))$ is continuous in $(t_i,t_{i+1})$. Given a Lagrangian for $\hh$, we  define the corresponding hybrid cost functional $C$ by
$$
C(\tau,q(\cdot),x(\cdot))=\sum_{i=0}^{M-1}\int_{t_{i}}^{t_{i+1}}L_{q_i}(t, x_i(t))\,dt .
$$
Let $(q_0,x_0)\in Q\times X_{q_0}$ be fixed.
We consider the hybrid optimization problem
\begin{equation*}
\begin{cases}
&\displaystyle{\min C(\tau,q(\cdot),x(\cdot))}, \\
&(\tau,q(\cdot),x(\cdot)) \textrm{ trajectory of }{\cal H},\\
&q(0)=q_0,~~x(0)=x_{0} .
\end{cases}
\eqno{(\mathrm{HP})}
\end{equation*}
Let $Q=\{q_1,\dots, q_k\}$. We define $h:Q\to \{1,\dots, k\}$ by $h(q_i)=i$.
For every $\varepsilon\geq 0$, we consider the optimization problem
\begin{equation*}
\begin{cases}
&\displaystyle{\min C(\tau,q(\cdot),x(\cdot))} +\varepsilon\, \mathrm{TV}(h\circ q(\cdot)), \\
&(\tau,q(\cdot),x(\cdot)) \textrm{ trajectory of }{\cal H},\\
&q(0)=q_0,~~x(0)=x_{0}.
\end{cases}
\eqno{(\mathrm{HP})_\varepsilon}
\end{equation*}

Casting Theorem~\ref{thmtv} in the language of hybrid systems, we obtain the following result.

\begin{theorem}\label{thmtvzeno}
Let $\hh$ be a hybrid system such that $X_q$ is a compact submanifold for every $q\in Q$, and such that the sets $G(q,q'), R((q,q'),x)$ are compact for every $((q,q'),x)\in E\times  X_q$, with $q\in Q$. 
Let $L$ be a Lagrangian for $\hh$ with corresponding cost functional $C$.
Assume that $(\tau^*,q^*(\cdot),x^*(\cdot))$ is a Zeno trajectory, optimal solution of $(\mathrm{HP})$.
For every $\varepsilon> 0$, the problem $(\mathrm{HP})_\varepsilon$ has at least one solution.
Moreover, for any solution $(\tau^\varepsilon, q^\varepsilon(\cdot), x^\varepsilon(\cdot))$ of $(\mathrm{HP})_\varepsilon$, we have
\begin{equation}\label{costizeno}
\lim_{\varepsilon\to 0} C(\tau^\varepsilon,q^\varepsilon(\cdot), x^\varepsilon(\cdot))=C(\tau^*,q^*(\cdot), x^*(\cdot)).
\end{equation}
\end{theorem}

The compactness assumption  on $X_q$ can be slightly weakened, and replaced by compactness of trajectories in each location.
The main idea of the proof is to  interpret  the role of the discrete part of the hybrid system in \eqref{csh} as a control. 
Since there are no final conditions, the proof is simplified with respect to the ones of Theorem~\ref{thmtv} and of Theorem~\ref{thm:rate}.

The rate of convergence in~(\ref{costizeno}) can be determined in the case where the rate of convergence of the switching times along the Zeno trajectory is known. We refer to Remark~\ref{rk:9292} in Section \ref{proofzeno} (end of the proof of Theorem \ref{thmtvzeno}) for a precise statement.

\begin{remark} 
In the definition of the hybrid system, one may now add a control. We do not provide the details. For such hybrid optimal control problems, assuming moreover that, in each location, the epigraph of extended velocities (defined by \eqref{defV}) is convex, that ${\bf U}$ is compact and that \eqref{hypub} holds true, the conclusion of Theorem~\ref{thmtvzeno} still holds true. In other words we have exactly the conclusion of Theorem \ref{thmtv} in the hybrid framework, including the convergence of trajectories.
\end{remark}

\begin{remark}
The problem of finding necessary and sufficient conditions for the existence of Zeno trajectories of a hybrid system has been firstly addressed in~\cite{zeno-numpb}, in which the authors dealt with the regularization of two specific hybrid systems: water tank and bouncing ball. Exploiting the specific geometry of the system, they introduced a family of regularized problems whose solution is ``close to" the Zeno trajectory. Their idea was either to introduce an additional variable whose role is to delay of $\varepsilon$ the time at which a switch takes place, or to introduce a spatial hysteresis. We refer to~\cite{ZJLS00,ZJLS01}
for a large number of examples of Zeno hybrid systems from the areas of modelling, simulation, verification, and control as well as for a list of references on the subject. 
We also refer to~\cite{AAS07,HLMR05,SJSL00} where conditions for the existence of Zeno solutions have been established.
The Zeno phenomenon for hybrid systems is related to so-called Zeno equilibria, which are invariant under the discrete (but not under the continuous) dynamics.
See also~\cite{LA13} for asymptotic stability of Zeno equilibra.
\end{remark}

\section{Proofs}\label{sec:proof}
\subsection{Proof of Theorem~\ref{thmtv}}
Before going into technical details, let us outline the proof of Theorem~\ref{thmtv}. First, the local controllability assumption implies the existence of an optimal solution $(u_\varepsilon,x_\varepsilon)$ of $\ocpe$ for any $\varepsilon\geq 0$. Second, thanks to the assumptions on the extended velocity sets and on the equiboundedness of trajectories, there exists an admissible control $w$ and a positive measurable function $\g$ such that the family $t\mapsto (f(x_\varepsilon(t),u_\varepsilon(t)),L(t,x_\varepsilon(t),u_\varepsilon(t)))$ converges to $t\mapsto(f(x_w(t),w(t)),L(t,x_w(t),w(t))+\g(t))$ for the weak star topology of $L^\infty$. Third, we use the optimality of $u_\varepsilon$ to prove that $w$ is optimal for $\ocp$. Fourth, we establish that $\g=0$, which implies that the Lagrangian cost along $u_\varepsilon$ converges to the Lagrangian cost at $u^*$. 
This fact is proved thanks to Lemma~\ref{dens} which exhibits a sequence of admissible controls $v_n$ for which $\mathrm{TV}(v_n)<+\infty$ and whose Lagrangian costs converge to the cost of $u^*$. To construct $v_n$, we use a topological result (Lemma~\ref{lemma-top-aux}), providing admissible controls steering any point of a neighborhood of the origin to $0$, with controls having bounded total variation.

We start by presenting the two auxiliary lemmas mentioned above, and then we proceed to the proof of the theorem. Note that the two lemmas do not require the convexity assumption of \eqref{defV} nor the \emph{a priori} estimate \eqref{hypub} on trajectories.

\begin{lemma}\label{lemma-top-aux}
Assume that $\Lie_0{\cal F}=\R^N$ and that the control system \eqref{dinamica} is small-time locally controllable at $0$. Then, there exists a neighborhood ${\cal N}$ of $0$ such that, for every $y\in {\cal N}$, there exists a piecewise constant control  $w_{y}:[0,\tau_{y}]\to {\bf U}$ steering \eqref{dinamica} from $y$ to $0$ in time $\tau_{y}$, with $\lim \tau_{y} =0$ as $y\to 0$.
\end{lemma}

\begin{proof}
Since \eqref{dinamica} is STLC at $0$, by \cite[Theorem~5.3 a-d]{grasse} we have that the
reversed control system
\begin{equation}\label{mdinamica}
\dot x=-f(x,u),\quad u\in{\cal U}, \tag{$-\Sigma$}
\end{equation}
associated with the dynamics $-f$
is also STLC at $0$.
As a consequence the time optimal map $\bar \Upsilon$ associated with system \eqref{mdinamica}, namely
$x_0\mapsto \bar \Upsilon(x_0)= \inf\{t>0\mid \dot x = - f(x,u), x(0)=0, x(t)=x_0\}$ is continuous at $0$ (see~\cite[Theorem 2.2]{stefani-holder}).
We denote by $\mathcal{A}(x, [0,T), -f)$, respectively by $\mathcal{A}(x, [0,T], -f)$, the set of points accessible from $x$ in time $t< T$, respectively  $t\leq T$, by trajectories of the control system \eqref{mdinamica}. We set $\mathcal{A}(x, -f)=\cup_{T>0}\mathcal{A}(x, [0,T), -f)$.
By definition of STLC there exists a neighborhood $\mathcal{N}$ of $0$ such that 
$\mathcal{N} \subset \mathcal{A}(0,-f)$.
Let $y \in \mathcal{N}$. By definition of the time optimal map, we have that
$
y \in \mathcal{A}(0,[0,\bar\Upsilon(y)], -f)\subset\mathcal{A}(0,[0,2\bar\Upsilon(y)), -f).
$
In particular this implies (see~\cite[Theorem~5.5]{grasse}) that
$y$ is \emph{normally reachable} (see~\cite[Definition~3.6]{grasse}) from $0$ in time less than $2\bar\Upsilon(y)$ for the control system \eqref{mdinamica}. Namely, there exist $q=q(y)\in\N$, $u_1,\dots u_q\in{\bf U}$ and positive numbers $t_1,\dots, t_q$ with $t_1+\dots + t_q < 2\bar\Upsilon(y)$, such that  
$y = \exp(-t_q f(\cdot, u_q))\circ\cdots\circ\exp(-t_1 f(\cdot, u_1)) (0)$.
Here, $\exp(tV)$ designates the flow at time $t$ of the vector field $V$.
Since $f$ is autonomous, we obtain
$\exp(t_1 f(\cdot, u_1))\circ\cdots\circ\exp(t_q f(\cdot, u_q)) (y)=0$.
Setting $\tau_y=t_1+\dots+t_q$ and defining $w_y:[0, \tau_y] : \to \R$ by
$$
w_y(t)=
\begin{cases}
u_1,&t\in[0,t_1],\\
u_2, &t\in[t_1,t_1+t_2],\\
\vdots & \\
u_q, &t\in[t_1+\dots+t_{q-1},\tau_y] ,
\end{cases}
$$
the lemma follows.
\end{proof}

 \begin{lemma}\label{dens}
Let $u:[0,t(u)]\to {\bf U}$  be a measurable control steering $x^0$ to $0$. Then there exists a countable family of controls $u_n:[0,t(u_n)]\to{\bf U}$ such that $\mathrm{TV}(u_n)<+\infty$ for every $n\in\N$, $u_n$ steers the control system \eqref{dinamica} from $x^0$ to $0$ in time $t(u_n)$ and
\begin{equation*}
\lim_{n\to +\infty}\|u_n-u\|_{L^1}=0.
\end{equation*}
\end{lemma}
Here, the $L^1$ norm is on $[0,+\infty)$, by extending $u$ (resp., $u_n$) by $0$ for $t>t(u)$ (resp., $t>t(u_n)$). Recall that $f(0,0)=0$, and thus this extension does not have any impact on admissible trajectories.

\begin{proof}
Consider a sequence of functions $v_n:[0,t(u)]\to{\bf U}$ with $\mathrm{TV}(v_n)<+\infty$ for every $n\in\N$ converging to $u$ in $L^1([0,t(u)],\R^m)$ for the strong topology and consider  the associated solutions $y_n(\cdot)$ of the Cauchy problem $\dot y_n =f(y_n, v_n)$, $y_n(0)=x^0$. Then the sequence  $y_n(\cdot)$ converges uniformly to the trajectory  $x_u(\cdot)$ associated with the control $u$ (see for instance~\cite[Theorem 3.4.1]{benlibro}).
In particular, $y_n(t(u))$ converge to $0$ as $n$ tends to $+\infty$. By Lemma~\ref{lemma-top-aux}, for $n$ sufficiently large, there exists a control $w_n:[0,\tau_n]\to {\bf U}$ which is piecewise constant, of bounded variation, steering $y_n(t(u))$ to $0$ in time $\tau_n$ and such that $\tau_n \to0$ as $n \to \infty$.
Define 
\bqnn
u_n(t)=\begin{cases}
v_n(t),&t\in[0,t(u)),\\
w_n(t-t(u)),&t\in[t(u),t(u)+ \tau_n),\\
0,&t>t(u)+\tau_n.
\end{cases}
\eqnn
By construction, $u_n$ steers $x^0$ to $0$ in time $t(u_n)=t(u)+\tau_n$ and, for every $n$, one has $\mathrm{TV}(u_n)<+\infty$. We extend $u$ to $[0,+\infty)$ by setting $u(t)=0$ for $t>t(u)$. 
Then
\begin{align*}
\int_0^{+\infty} & |u_n(s)-u(s)|\,ds \\
& = \int_0^{t(u)} |v_n(s)-u(s)|\, ds+\int_{t(u)}^{t(u)+\tau_n} |w_n(s-T)|\, ds \\
&\leq \int_0^{t(u)} |v_n(s)-u(s)|\,ds+ \tau_n\max_{z\in{\bf U}}|z|,
\end{align*}
which converges to zero since $\tau_n \to 0$ as $n\to\infty$ and for to the strong convergence of $v_n$ towards $u$ in $L^1$.
\end{proof}

Let us now prove Theorem~\ref{thmtv}. The proof follows the lines of~\cite[Theorem 5.14 and 6.15]{trelat_book}.

\begin{proof}[Proof of Theorem~\ref{thmtv}]
First of all, by Lemma \ref{lemma-top-aux}, there exists a control $u:[0,t(u)]\to {\bf U}$ steering $x_0$ to $0$ and having bounded variation. Therefore, the existence of an optimal solution of $\ocpe$ follows from Theorem~\ref{exth} in Appendix \ref{sec:app}. 

Let $x_\varepsilon(\cdot)$ be any optimal solution of $\ocpe$, associated with a control $u_\varepsilon:[0,t(u_\varepsilon)]\to {\bf U}$.
Set $\tilde x_\varepsilon(t) = (x_\varepsilon(t), \int_0^t L(s,x_\varepsilon(s),u_\varepsilon(s)) ds)$. Then the triple $(\tilde x_\varepsilon,u_\varepsilon, \gamma_\varepsilon)$ with $\gamma_\varepsilon \equiv 0$ is a solution of 
$$
\min_{u\in \mathcal{U}, \gamma \geq 0} \left(\int_0^{t(u)} (L(s,x,u) + \gamma(s) )ds + \varepsilon  \mathrm{TV}(u) \right) 
$$
subject to 
\begin{equation}\label{eq:eeee}
  \dot x = f( x,u),\quad \dot x_{N+1}=L(t, x,u) + \gamma,
\end{equation}
with initial conditions $x(0)=x_0$, $x_{N+1}(0)=0$, and final conditions $x(t(u))=0$, $x_{N+1}(t(u))\geq 0$.
Denote by $\tilde f(t,x,u,\gamma)=(f(x,u),L(t, x,u)+\gamma)$ the augmented dynamics of~\eqref{eq:eeee} which are convex by Assumption~\eqref{defV}.

Thanks to Assumption \eqref{hypub}, the sequence $t(u_\varepsilon)$ is bounded and converges, up to some subsequence, to $t_1>0$ as $\varepsilon$ tends to $0$. Hence, given $\delta >0$ there exists $\varepsilon_0>0$ such that $|t(u_\varepsilon)-t_1|<\delta$ for every $\varepsilon\in[0,\varepsilon_0]$ in the chosen subsequence. 
Since $f(0,0)= 0$, we extend  $x_\varepsilon$ and $u_\varepsilon$ to $[t(u_\varepsilon), t_1+\delta]$ by $0$.
By Assumption \eqref{hypub}, the trajectories $x_\varepsilon(\cdot)$ are uniformly bounded, and hence the family of functions $s\mapsto \tilde f(s, x_\varepsilon(s),u_\varepsilon(s),0)$ is bounded in $L^\infty([0,t_1+\delta],\R^{N+1})$. Thus, up to some subsequence, it converges to some function $g\in L^\infty([0,t_1+\delta],\R^{N+1})$ for the weak star topology. We define
$$
\tilde x(t)= \tilde x_0+\int_0^tg(s)\,ds,\quad\tilde x_0=(x_0,0).
$$
By construction, $t\mapsto \tilde x(t)$ is absolutely continuous. Moreover, the family $\tilde x_\varepsilon(t)$ converges uniformly to $\tilde x(\cdot)$ on $[0,t_1+\delta]$. By the convexity assumption~\eqref{defV} the absolutely continuous function $\tilde x(\cdot)$ is also a trajectory of~\eqref{eq:eeee} (see, for instance~\cite[Corollary 3.3.2]{benlibro}), in 
particular
there exists an admissible control $w:[0,t(w)]\to{\bf U}$ and a positive measurable function $\gamma:[0,t(w)]\to \R$ such that 
$t \mapsto \tilde x(t) := (x_w(t), \int_0^tL(s, x_w(s),u_w(s))+ \gamma(s) ds )$ is the associated solution of~\eqref{eq:eeee}.

It remains to prove that $x_w(\cdot)$ is optimal for $\ocp$.
For every admissible control $v\in{\cal U}$ satisfying $\mathrm{TV}(v)<+\infty$, we have (note that $\gamma(\cdot)\geq 0$)
\begin{align}
& \int_0^{t(w)+\delta} L(t,x_w(t),w(t))\,dt \nonumber\\
&  \leq\int_0^{t(w)+\delta}(L(t,x_w(t),w(t))+\g(t))\,dt\nonumber\\
&\leq \limsup_{\varepsilon\to 0}\left(\int_0^{t(w)+\delta}L(t,x_\varepsilon(t),u_\varepsilon(t))+ \varepsilon\, \mathrm{TV}(u_\varepsilon)\right)\label{eq:12312}\\
&\leq\int_0^{t(v)}L(t,x_v,v)\,dt \nonumber \\
&\quad +\int_{t(w)}^{t(w)+\delta}(L(t,x_w(t),w(t))+\g(t))\,dt.
\nonumber
\end{align}
Hence, for every $\delta>0$ and every admissible $v$ as above, we have
\begin{align*}
\int_0^{t(w)} & L(t,x_w(t),w(t))\,dt\\
& \leq \int_0^{t(v)}L(t,x_v,v)\,dt +\int_{t(w)}^{t(w)+\delta}\g(t)\,dt .
\end{align*}
Since $\delta>0$ was taken arbitrary, we conclude that
$$
\int_0^{t(w)}L(t,x_w(t),w(t))\,dt\leq \int_0^{t(v)}L(t,x_v,v)\,dt ,
$$
for every admissible control $v\in{\cal U}$ satisfying $\mathrm{TV}(v)<+\infty$
Using Lemma~\ref{dens} and the dominated convergence theorem, we infer that the inequality above holds true as well for any possible admissible control $v\in{\cal U}$ (not necessarily of bounded variation). Therefore, $w$ is the optimal control solution of $\ocp$. 

Finally, to prove \eqref{convcost}, it suffices to show that $\g=0$. By optimality of $u_\varepsilon$, we have
\begin{align*}
\int_0^{t(u_\varepsilon)} & L(s,x_\varepsilon(s),u_\varepsilon(s))\,ds\\ 
& \leq \int_0^{t(u_\varepsilon)}L(s,x_\varepsilon(s),u_\varepsilon(s))\, ds+\varepsilon\, \mathrm{TV}(u_\varepsilon)\\
& \leq \int_0^{t(v)}L(s,x_v(s),v(s)))\,ds+\varepsilon\, \mathrm{TV}(v) ,
\end{align*}
for any admissible control $v$ such that $\mathrm{TV}(v)<+\infty$. 
Letting $\varepsilon$ tend to $0$, we deduce that
\begin{align*}
\int_0^{t(w)} &  ( L(s,x_w(s),w(s))+\g(s) ) \,ds\\
& \leq   \int_0^{t(v)}L(s,x_v(s),v(s)))\,ds.
\end{align*}
Finally, since $w$ is optimal for $\ocp$, we conclude that $\g=0$.
\end{proof}

\subsection{Proof of Theorem~\ref{thm:rate}}
We start with the following lemma.

\begin{lemma}\label{mainth}
Assume that the control system \eqref{dinamica} satisfies \condition\ at $0$.
Then, for every $\eta>0$ sufficiently small, there exists an admissible control $v_\eta:[0,t(v_\eta)]\to{\bf U}$ satisfying $\mathrm{TV}(v_\eta)<+\infty$, whose corresponding trajectory is denoted by $x_\eta(\cdot)$, such that 
\begin{equation}\label{costconv}
\lim_{\eta\to 0}\int_0^{t(v_\eta)}L(t,x_\eta(t),v_\eta(t))\,dt = \int_0^{t(u^*)}L(t,x^*(t),u^*(t))\,dt ,
\end{equation}
and 
\begin{align*}
\lim_{\eta\to 0}|t(v_\eta) -t(u^*)|&=\lim_{\eta\to 0} \|v_\eta -  u^*\|_{L^{1}}\\
&=\lim_{\eta\to 0} \|x_\eta(\cdot)-x^*(\cdot)\|_{\infty} = 0 .
\end{align*}
Moreover, under the additional assumption that the time-optimal map is $\con^{0,\alpha}$ for some $\alpha \in (0,1]$ in a neighborhood of $0$, there exists $C>0$ such that
\begin{align}
\int_0^{t(v_\eta)}L(t,x_\eta(t),v_\eta(t))\,dt- & \int_0^{t(u^*)}L(t,x^*(t),u^*(t))\,dt \nonumber\\
&\leq C\eta^{\alpha}.\label{costconveta}
\end{align}
\end{lemma}

\begin{proof}
Let ${\cal N}$ be the neighborhood of $0$ in $\R^N$ and $M$ be the constant given by Definition~\ref{def:omega}. Without loss of generality we can assume that  ${\cal N}$ is bounded.
Fix $\eta_{0}$  such that $x^*(s)\in {\cal N}$, for every $s\geq t(u^*)-\eta_0$.
By condition \condition, there exists a control $w_{\eta}$ steering $x^{*}(t(u^*)-{\eta})$ to $0$ in time $\tau_{\eta} \leq M \Upsilon(x^*(T-\eta))$  
with $\mathrm{TV}(w_{\eta}) \leq M$. We define $v_{\eta}$ by
\begin{align}
&v_{\eta}(t) =\nonumber \\
&=
\begin{cases}
u^{*}(t) &\mbox{ for } t \in [0,t(u^*)-{\eta}),\\
w_{\eta}(t-t(u^*)+{\eta}) & \mbox{ for } t \in [t(u^*)-{\eta},t(u^*)-{\eta}+\tau_\eta),\\
0, &\mbox{ for } t>t(u^*)-{\eta}+\tau_\eta ,
\end{cases}\label{eq:veta}
\end{align}
and let $x_{\eta}(\cdot)$ be the corresponding trajectory, starting from $x_0$ (see Figure~\ref{fig:lemma2}).
\begin{figure}[h]
\begin{center}
\includegraphics[width=13cm]{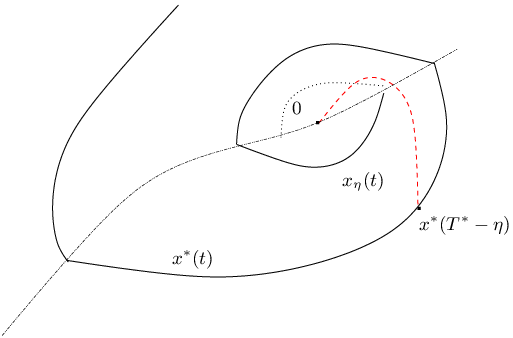} 
\caption{The trajectory $x_\eta(\cdot)$ associated with the control $v_\eta$.}
\label{fig:lemma2}
\end{center}
\end{figure}
By construction, we have $\mathrm{TV}(v_{\eta}) \leq \mathrm{TV}(u^{*}|_{[0,t(u^*)-{\eta}]}) + M$.
If $\mathrm{TV}(u^*)<+\infty$ or $u^*$ is chattering in the sense of Definition~\ref{def:chat}, then $\mathrm{TV}(v_\eta)<+\infty$. 
We have $\tau_\eta\to 0$ as $\eta \to 0$, since $\Upsilon$ is upper semi-continuous. Hence $v_\eta\to u^*$ almost everywhere, and for some subsequence, we have
$$
\lim_{\eta\to 0} T_{\eta} = t(u^*) \quad \mbox{ and } \quad \lim_{\eta\to 0} \|v_{\eta} - u^{*}\|_{L^{1}} =0.
$$
Now, set $\mathcal{X}_{0} = \{x^{*}(t)\ \mid\ t\in[0,t(u^*)-\eta_0]\} \cup \overline{{\cal N}}$,  
$C_{1} = \sup_{\mathcal{X}_{0} \times \mathbf{U}} |\partial_{x}f|$, and    
$C_{2} = \sup_{\mathcal{X}_{0} \times \mathbf{U}}  |\partial_{u}f|$.  For every $t\geq 0$ and for every $\eta\in (0,\eta_{0})$, we have
\begin{align*}
|x_{\eta}(t) & - x^{*}(t)| \\
& = \left| \int_{0}^{t} f(x_{\eta}(s),v_{\eta}(s))\, ds - \int_{0}^{t} f(x^{*}(s),u^{*}(s))\, ds   \right|\\
&\leq \int_{0}^{t}|f(x_{\eta}(s),v_{\eta}(s)) -  f(x^{*}(s),v_\eta(s))| \, ds\\
&\quad +\int_0^t|f(x^{*}(s),v_\eta(s))-f(x^{*}(s),u^{*}(s))| \, ds   \\
&\leq  C_{1} \int_{0}^{t} |x_{\eta}(s) - x^{*}(s)| \, ds + C_{2}\|u^{*}-v_{\eta}\|_{L^{1}},
\end{align*}
and thus, by the Gronwall lemma, we get that $\|x_{\eta}(\cdot) - x^{*}(\cdot)\|_{\infty} \leq  C_{2} \|u^{*}-v_{\eta}\|_{L^{1}} e^{C_{1}\bar T}$.
In particular $\lim_{\eta \to 0} \|x_{\eta}(\cdot) - x^{*}(\cdot)\|_{\infty} =0$.

Finally, let us prove \eqref{costconv}. By continuity of $L$,  there exist constants $c\in\R$ and $\bar C>0$ such that $L(t,x^*(t),u^*(t))\geq c$ for almost every $t\in [0,t(u^*)]$, and $|L(t,x,u)| \leq \bar C$ for almost every $(t,x,u)\in [0,\bar T]\times \mathcal{X}_{0} \times \mathbf{U}$.
Then, we have
\begin{eqnarray*}
0 & \leq & \int_{0}^{t(u^*)-\eta+\tau_\eta}L(t,x_\eta(t), v_\eta(t))\,dt \\
& &-\int_{0}^{t(u^*)}L(t,x^*(t),u^*(t))\,dt \\
&=& \int_{t(u^*)-\eta}^{t(u^*)-\eta+\tau_\eta}L(t,x_\eta(t), v_\eta(t))\,dt\\
& & -\int_{t(u^*)-\eta}^{t(u^*)}L(t,x^*(t),u^*(t))\,dt \\
&\leq& \bar C \tau_{\eta}   -  c\eta ,
\end{eqnarray*}
which implies  \eqref{costconv}. To prove~\eqref{costconveta} it suffices to note that $\tau_{\eta} \leq M\Upsilon(x^{*}(t(u^*)-\eta)) \leq C \eta^{\alpha}$.
\end{proof}

We are now in a position to prove Theorem~\ref{thm:rate}.

\begin{proof}[Proof of Theorem~\ref{thm:rate}] 
Assumption \condition\ implies in particular the existence of bounded variation controls steering the control system \eqref{dinamica} from any initial condition in the neighborhood $\mathcal{N}$ to the origin. Hence, from Theorem~\ref{exth} in Appendix \ref{sec:app}, the problem $\ocpe$ has at least one solution.

Let $x_{\varepsilon}(\cdot)$ be an arbitrary solution of $\ocpe$, associated with a control $u_{\varepsilon}:[0,T_{\varepsilon}] \to {\bf U}$.
Let $n_{0}\in\N$ be such that $x^*(t_n) \in \mathcal{N}$ for every $n \geq n_{0}$.
We apply Lemma~\ref{mainth} with $\eta = t(u^*) - t_{n}$ and we denote, for simplicity, $u_{n}$ the control $v_{t(u^*)-t_n}$.
Note that  $\mathrm{TV}(u_n) \leq n+M$.
By optimality of $u_{\varepsilon}$ for $\ocpe$, we have
\begin{align*}
 \int_0^{T_\varepsilon} &  L(t,x_\varepsilon(t),u_\varepsilon(t))\,dt \\
 &\leq \int_0^{T_\varepsilon} L(t,x_\varepsilon(t),u_\varepsilon(t))\,dt + \varepsilon\, \mathrm{TV}(u_\varepsilon)\\
 &\leq \int_0^{t_n+\tau_n}L(t,x_n(t),u_n(t))\,dt + \varepsilon\, \mathrm{TV}(u_n)\\
 & \leq \int_0^{t(u^*)}L(t,x^*(t),u^*(t))\,dt \\
 &\quad + C |t(u^*)-t_n|^\alpha + \varepsilon (n+M).
\end{align*}
Now, by Assumption $(ii)$ made in the statement of the Theorem, we have $|t(u^*)-t_n|^\alpha = \mathrm{O}(n^{-\alpha\beta})$, and choosing $n = \mathrm{O}(\varepsilon^{-\frac{1}{1+\alpha\beta}})$, we infer that 
\begin{align*}
\int_0^{T_\varepsilon} &  L(t,x_\varepsilon(t),u_\varepsilon(t))\,dt  - \int_0^{t(u^*)}L(t,x^*(t),u^*(t))\,dt \\  
& \leq C |t(u^*)-t_n|^\alpha + \varepsilon (n+M) = \mathrm{O}\left(\varepsilon^{\frac{\alpha\beta}{1+\alpha\beta}}\right).
\end{align*}
This concludes the proof.
\end{proof}

\subsection{Proof of Theorem~\ref{robbinseps}}\label{proofr}
We start with the following existence result.

\begin{lemma}\label{existence} 
Given any $\varepsilon>0$, the problem $\ocpse$ has at least one solution.
\end{lemma}

\begin{proof}
Let $\varepsilon>0$ be fixed. First of all, remark that if there exists no admissible trajectory such that $\mathrm{TV}(X_u)<+\infty$, then the functional $u\mapsto \int_0^TL(s,x(s),u(s))\,ds+\varepsilon\, \mathrm{TV}(X_u)$ is infinite and there is nothing to prove.
Otherwise, let $I<+\infty$ denote the infimum in $\ocpse$.
We consider a minimizing sequence of admissible controls $u_n:[0,t(u_n)]\to{\bf U}$, with corresponding trajectories denoted by $x_n(\cdot)$, such that
$$
\lim_{n\to\infty}\left(\int_0^{t(u_n)}L(s,x_n(s),u_n(s))\,ds+ \varepsilon\, \mathrm{TV}(X_{u_n})\right)=I.
$$ 
Since the sequence $(t(u_n))_{n\in\N}$ is bounded by Assumption \eqref{hypub}, we can assume that $t(u_n)$ converges (up to some subsequence) to some $t_\varepsilon>0$. Using  $f(0,0)=0$ we extend $u_n$ to $[t(u_n),t_\varepsilon+\delta]$ by $0$ for $\delta>0$.
Reasoning as in the proof of Theorem~\ref{thmtv}, up to some subsequence, there exist a positive measurable function $\g:[0,t_\varepsilon+\delta]\to\R$ and a measurable control $w:[0, t_\varepsilon+\delta]\in{\bf U}$, with corresponding trajectory $x_w(\cdot)$, such that
$x_n(\cdot)$ converges to $x_w(\cdot)$ uniformly on $[0, t_\varepsilon+\delta]$ and $L(\cdot, x_n(\cdot),u_n(\cdot))$ converges to $L(\cdot, x_w(\cdot), w(\cdot))+\g(\cdot)$ in $L^\infty(0,t_\varepsilon+\delta)$ for the weak star topology.
By uniform convergence of trajectories, $w:[0,t_\varepsilon]\to{\bf U}$ is admissible, that is, $x_w(t_\varepsilon)=0$ and $x_w(t)\in\con$ for every $t$.
Up to some subsequence, by dominated convergence, we can assume that $X_{u_n}(\cdot)={\bf 1}_{\partial \con}(x_n(\cdot))$ converges to $X_w(\cdot)={\bf 1}_{\partial \con}(x_w(\cdot))$ in $L^1(0,t_\varepsilon+\delta)$. Moreover, since $x_n(t)=0$ on $[t(u_n),t_\varepsilon+\delta]$, we have $\mathrm{TV}(X_{u_n}|_{[0,t_\varepsilon+\delta]})= \mathrm{TV}(X_{u_n}|_{[0,t(u_n)]})$. Therefore,
\begin{align*}
\int_{0}^{t_\varepsilon+\delta}& L(t,x_n(t),u_n(t))\,dt+\varepsilon\, \mathrm{TV}(X_{u_n}|_{[0,t_\varepsilon+\delta]})\\
& = \int_{0}^{t(u_n)}L(t,x_n(t),u_n(t))\,dt+\varepsilon\, \mathrm{TV}(X_{u_n}|_{[0,t(u_n)]})\\
& \quad +\int_{t(u_n)}^{t_\varepsilon+\delta}L(t,x_n(t),u_n(t))\,dt.
\end{align*}
We infer that
\begin{align*}
\limsup_{n\to\infty}& \left(\int_{0}^{t_\varepsilon+\delta}L(t,x_n(t),u_n(t))\,dt+\varepsilon\, \mathrm{TV}(X_{u_n}|_{[0,t_\varepsilon+\delta]})\right)\\
& \leq  I+\limsup_{n\to\infty}\int_{t(u_n)}^{t_\varepsilon+\delta}L(t,x_n(t),u_n(t))\,dt\\
&=I+\int_{t_\varepsilon}^{t_\varepsilon+\delta}(L(t,x_w(t),w(t))+\g(t))\,dt.
\end{align*}
Besides, by lower semicontinuity of $\mathrm{TV}(\cdot)$, we have $\mathrm{TV}(X_w)<+\infty$ and
\begin{align*}
& \liminf_{n\to\infty} \left(\int_{0}^{t_\varepsilon+\delta}L(t,x_n(t),u_n(t))\,dt  +\varepsilon\, \mathrm{TV}(X_n|_{[0,t_\varepsilon+\delta]})\right)\\
& \quad \geq \liminf_{n\to\infty}\int_{0}^{t_\varepsilon+\delta}L(t,x_n(t),u_n(t))\,dt + \varepsilon\, \mathrm{TV}(X_w|_{[0,t_\varepsilon+\delta]}) \\
&  \quad \geq \int_0^{t_\varepsilon+\delta}(L(t,x_w(t),w(t))+\g(t))\,dt+\varepsilon\, \mathrm{TV}(X_w|_{[0,t_\varepsilon]}).
\end{align*}
Finally, we obtain that $I\geq \int_0^{t_\varepsilon}(L(t,x_w(t),w(t))+\g(t))\,dt+\varepsilon\, \mathrm{TV}(X_w|_{[0,t_\varepsilon]})$.
Since $w$ is admissible, we have $x_w(t_\varepsilon)=0$ and there holds
$I\leq  \int_0^{t_\varepsilon}L(t,x_w(t),w(t))\,dt + \varepsilon\, \mathrm{TV}(X_w|_{[0,t_\varepsilon]})$.
Therefore, since $\gamma \geq 0$, we infer that $\int_0^{t_\varepsilon}\g(t)\,dt=0$ and $w:[0,t_\varepsilon]\to{\bf U}$ is optimal for $\ocpse$.
\end{proof}

\begin{lemma}\label{densss}
Assume condition (i) of Theorem~\ref{robbinseps}. 
Let $\bar u:[0,\bar t]\to{\bf U}$ be an admissible control 
such that the corresponding trajectory $\bar x(\cdot )$ satisfies $\bar x(t)\in\con$ for every $t$ and 
$\{t\mid \bar x(t)\in\partial\con\}=\{\bar t,t_1,t_2,\dots\}$ with $\lim_{n\to\infty}t_n=\bar t$. Then, there exists a sequence $(u_k)_{k\in\N}$ of ${\cal U}$ of admissible controls,
such that $u_k$ converges to $\bar u$ in $L^1$, the corresponding trajectories $x_k(\cdot)$ satisfy $x_k(t)\in\con$ for every $t$, and $\mathrm{TV}(X_{u_k})<+\infty$. 
\end{lemma}

\begin{proof}
Fix $k > 0$. Recall that condition $(i)$ states that there exists a neighborhood $\mathcal{N}$ of $0$ such that
$\mathcal{N}\cap \mathring{\con} \subset \mathcal{A}^\mathcal{C}(0,(0,1/k),-f)$.
 By assumption,
for almost every $\eta>0$, the point $x^\eta=\bar x(\bar t-\eta)$ belongs to the interior of $\con$. Hence for almost every $\eta>0$ sufficiently small we have $x^\eta\in{\cal N}\cap \mathring{\con}\subset {\cal A}^\con(0,(0,\eta),-f)$. Then, there exists a control $w^\eta:[0,\tau_\eta]\to {\bf U}$, with $\tau_\eta\leq \eta$, such that the solution $y(\cdot)$ of the Cauchy problem $\dot y=-f(y, w)$, $y(0)=0$, satisfies $y(t)\in \mathring{\con} $ for every $t\in(0,\tau_\eta]$ and $y(\tau_\eta)=x^\eta$.
Reversing time, since the dynamics is autonomous, we get that $z(\tau_\eta)=0$, where $z(\cdot)$ is the solution of the Cauchy problem $\dot z=f(z,w)$, $z(0)=x^\eta$, and $z(t)\in\mathring\con$ for every $t\in [0,\tau_\eta)$.
Let $T>0$. We extend the control $\bar u$ to $[\bar t, \bar t+T]$ by setting $\bar u= 0$. We define
$$
u_\eta(t)=
\begin{cases}
\bar u(t),& t\in[0,\bar t-\eta]\\
w^\eta(t-\bar t+\eta),& t\in [\bar t-\eta,\bar t-\eta+\tau_\eta]\\
0,&t>\bar t-\eta+\tau_\eta.
\end{cases}
$$
Since $\tau_\eta$ converges to $0$ as $\eta\to 0$, $u_\eta$ converges to $\bar u$ in $L^1(0,\bar t+T)$. Therefore, the sequence of corresponding trajectories $x_{\eta}(\cdot)$ converges uniformly to $\bar x(\cdot)$ on $[0,\bar t+T]$ (see Figure \ref{fig:ro}).
\begin{figure}[h]
\begin{center}
\includegraphics[width=13cm]{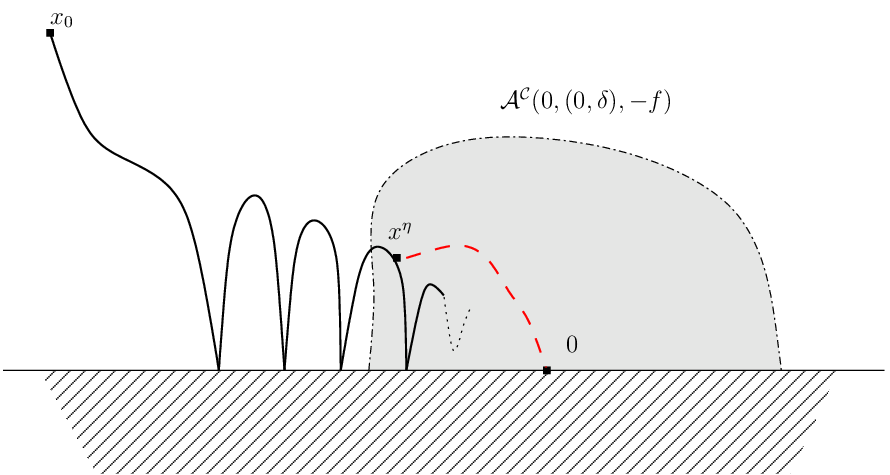} 
\caption{The trajectory $x_{\eta}(\cdot)$ associated with the control $u_\eta$.}
\label{fig:ro}
\end{center}
\end{figure}
Thus $L(\cdot,x_{\eta}(\cdot),u_\eta(\cdot))$ converges to $L(\cdot,\bar x(\cdot),\bar u(\cdot))$ strongly in $L^\infty(0,\bar t+T)$.
Set $T_\eta=\bar t-\eta+\tau_\eta$. By construction, we have $x_{\eta}(T_\eta)=0$ and $\mathrm{TV}({\bf 1}_{\partial\con}(x_{\eta})|_{[0,T_\eta]})=\mathrm{TV}({\bf 1}_{\partial\con}(\bar x )|_{[0,\bar t-\eta]})<+\infty$. Finally, the convergences above imply that $\int_0^{T_\eta}L(t,x_{\eta}(t),u_\eta(t))\,dt$ converges to $\int_0^{\bar u}L(t,\bar x(t),\bar u(t))\,dt$ as $\eta\to 0$. The statement follows by taking a sequence $\eta = 1/k$ for $k\in\N$ sufficiently large. 
\end{proof}

We are now in a position to prove Theorem~\ref{robbinseps}.

\begin{proof}[Proof of Theorem~\ref{robbinseps}]
Let $x_\varepsilon(\cdot)$ be any optimal solution of $\ocpse$, associated with a control $u_\varepsilon$ (existence is ensured by Lemma~\ref{existence}). We make the same reasoning as in the proof of Theorem~\ref{thmtv}. 
 
Let $t(u_\varepsilon)$ converge (up to some subsequence) to some $t_1>0$. Let $\delta>0$ be arbitrary. We extend $u_\varepsilon$ to $[0, t_1+\delta]$ by $0$. As in the previous proofs, there exists an admissible control $w:[0,t_1+\delta]\to{\bf U}$, with corresponding trajectory $x_w(\cdot)$, and a positive measurable function $\g:[0,t_1+\delta]\to\R$ such that $x_\varepsilon(\cdot)$ converges to $x_w(\cdot)$ uniformly on $[0,t_1+\delta]$, and $L(\cdot, x_\varepsilon(\cdot),u_\varepsilon(\cdot))$ converges to $L(\cdot, x_w(\cdot),w(\cdot))+\g(\cdot)$ in $L^\infty(0,t_1+\delta)$ for the weak star topology.
Replacing the total variation of controls $t\mapsto u(t)$ with the total variation of $t\mapsto X_u(t)={\bf 1}_{\partial\con}(x_u(t))$ in~\eqref{eq:12312}, we get that, for every admissible control $v:[0,t(v)]\to{\bf U}$ such that $\mathrm{TV}(X_v)<+\infty$, there holds
\begin{align*}
\int_0^{t_1} & L(t,x_w(t),w(t))\,dt\\
& \leq \int_0^{t(v)}L(t,x_v(t),v(t))\,dt+\int_{t_1}^{t_1+\delta}\g(t)\,dt.
\end{align*}
Since $\delta>0$ is arbitrary, letting $\delta$ tend to zero we conclude that, for every $v$ as above,
\begin{equation}\label{fghj}
\int_0^{t_1}L(t,x_w(t),w(t))\,dt\leq \int_0^{t(v)}L(t,x_v(t),v(t))\,dt.
\end{equation}
We apply Lemma~\ref{densss} to $\bar u=u^*$ and we denote by $u_k$ the corresponding sequence. Then, taking the inequality~\eqref{fghj} with $v=u_k$ and letting $k$ tend to $+\infty$, we obtain that $w$ is optimal for $\ocps$. In order to establish \eqref{convcostvinc}, it remains to prove that $\g|_{[0,t_1]}\equiv 0$. To this aim, let $v:[0,t(v)]\to{\bf U}$ be an admissible control such that $\mathrm{TV}(X_v)<+\infty$. Then, by optimality of $u_\varepsilon$ for $\ocpse$, we have
\begin{align*}
\int_0^{t_1+T} & (L(t,x_w(t),w(t))+\g(t))\,dt \\
& =\lim_{\varepsilon\to 0}\int_0^{t_1+T}L(t,x_\varepsilon(t),u_\varepsilon(t))\,dt \\ 
&\leq \limsup_{\varepsilon \to 0}\left(\int_0^{t_1+T}L(t,x_\varepsilon(t),u_\varepsilon(t))\,dt+\varepsilon\, \mathrm{TV}(X_{u_\varepsilon})\right)\\
&\leq \int_0^{t(v)}L(t,x_v(t),v(t))\,dt\\
& \quad+\limsup_{\varepsilon \to 0}\int_{t(u_\varepsilon)}^{t_1+T}L(t,x_\varepsilon(t),u_\varepsilon(t))\,dt\\
&= \int_0^{t(v)}L(t,x_v(t),v(t))\,dt \\
& \quad+ \int_{t_1}^{t_1+T}(L(t,x_w(t),w(t))+\g(t))\,dt,
\end{align*}
which gives $\int_0^{t_1}(L(t,x_w(t),w(t))+\g(t))\,dt\leq \int_{0}^{t(v)}L(t,x_v(t),v(t))\,dt$.
Again, let $(u_k)_{k\in\N}$ be the sequence provided by Lemma~\ref{densss} with $\bar u= u^*$. Then, since $w$ is optimal for $\ocps$, the inequality above with $v=u_k$ implies that $\int_0^{t_1}\g(t)\,dt=0$, which gives $\g_{[0,t_1]}=0$.
\end{proof}

\subsection{Proof of Theorem~\ref{thmtvzeno}}\label{proofzeno}
We first prove an auxiliary lemma.   

\begin{lemma}\label{lemzeno}
Let $\hh$ be a hybrid system
and let $L$ be a Lagrangian for $\hh$ with corresponding cost functional $C(\cdot)$.
Assume that $(\tau^*,q^*(\cdot),x^*(\cdot))$ is a Zeno solution of $\mathrm{(HP)}$.
Let $\tau^*=\{\tau^*_i\}_{i=0}^\infty$. Define the  sequence of trajectories $(\tau^n,q^n(\cdot), x^n(\cdot))$ by
\begin{itemize}
\item $\tau^n=\{\tau^*_0,\tau^*_1,\dots, \tau^*_n,\tau^*_\infty\}$;
\item $q^n(t)=q^*(t)$ for every $t\in[0,\tau^*_n)$, $q^n(t)\equiv q^*(\tau^*_n)$ for $t\in[\tau^*_n,\tau^*_\infty]$;
\item $x^n(t)=x^*(t)$ for every $t\in[0,\tau^*_n]$, and on $[\tau^*_n,\tau^*_\infty]$ the (continuous) trajectory $x^n(\cdot)$ is solution of $\dot x^n(t)=f_{q^*(\tau^*_n)}(x^n(t))$ almost everywhere.
\end{itemize} 
Then:
\begin{align}
&
\sup_{\tau^*_0 \leq t \leq \tau^*_n}||x^n(t)-x^*(t)||\leq \mathrm{O}(\tau^*_\infty-\tau^*_n) , \label{convx} \\
&C(\tau^n,q^n(\cdot), x^n(\cdot))-C(\tau^*,q^*(\cdot), x^*(\cdot))\leq \mathrm{O}(\tau^*_\infty-\tau^*_n). \label{convcosti}
\end{align}
\end{lemma}

\begin{proof}
Since $q^n(t)$ converges to $q^*(t)$ almost everywhere in $[0,\tau^*_\infty]$, by standard convergence results (see for instance~\cite[Theorem~1 p.~57]{sontag}) we deduce \eqref{convx}. For \eqref{convcosti}, note that, since the Lagrangian is continuous, there exist positive constants $\tilde c$ and $c$, satisfying $\tilde c-c>0$, such that $\int_{\tau^*_{n}}^{\tau_\infty^*}L_{q^*(\tau^*_{n})}(t,x^n(t))\,dt\leq \tilde c (\tau_\infty^*-\tau^*_{n})$ for every $n$, and $L_{q^*(\tau^*_i)}(t,x^*(t))\geq c$ almost everywhere in $[\tau^*_i,\tau^*_{i+1}]$ for every $i$.
Therefore,
\begin{align*}
0 &\leq C(\tau^n,q^n(\cdot), x^n(\cdot))-C(\tau^*,q^*(\cdot), x^*(\cdot))\\
&=\int_{\tau^*_{n}}^{\tau_\infty^*}L_{q^*(\tau_{n_i})}(t,x^n(t))\,dt-\sum_{i=n}^\infty\int_{\tau^*_i}^{\tau^*_{i+1}}L_{q^*(t)}(t,x^*(t))\,dt\\ 
&\leq \tilde c (\tau_\infty^*-\tau^*_{n})-c\sum_{i=n}^\infty(\tau^*_{i+1}-\tau^*_i)
=(\tilde c-c)(\tau_\infty^*-\tau^*_{n}).
\end{align*}
This concludes the proof of \eqref{convcosti}.
\end{proof}

Consider the set ${\cal T}_{M}$ of trajectories  of a hybrid system having at most $M$ switchings. If  $M<+\infty$ these trajectories are non-Zeno. We say that two non-Zeno trajectories  have the {\it same history} if they visit the same locations in the same sequence. Having the same history is an equivalent relation in ${\cal T}_{M}$ and the number of equivalent classes in ${\cal T}_{M}$ is finite. 

Let us now prove Theorem~\ref{thmtvzeno}.

\begin{proof}[Proof of Theorem~\ref{thmtvzeno}]
By compactness of the location $X_{q_0}$ there exists at least one trajectory starting at $(q_0,x_0)$ and having only a finite number of location switchings, for every $\varepsilon>0$ the functional to minimize in $\mathrm{(HP)}_\varepsilon$ is finite.  Therefore there exists $M_{\varepsilon}$ such that any solution of
\begin{equation*}
\begin{cases}
&\displaystyle{\min C(\tau,q(\cdot),x(\cdot))} +\varepsilon\, \mathrm{TV}(h\circ q(\cdot)), \\
&(\tau,q(\cdot),x(\cdot)) \in {\cal T}_{M_{\varepsilon}},\\
&q(0)=q_0,~~x(0)=x_{0},
\end{cases}
\eqno{(\mathrm{HP})^{'}_\varepsilon}
\end{equation*}
where the minimization runs over all possible trajectories having only a finite number of switchings $M_{\varepsilon}$
is also a solution of
$\mathrm{(HP)}_\varepsilon$.
Now consider a minimizing sequence for $(\mathrm{HP})^{'}_\varepsilon$. Then, up to some subsequence, we can assume that all trajectories have the same history. 
Hence the penalization term of total variation is constant along the chosen subsequence and the problem is then reduced to that of minimizing the Lagrangian cost $C(\cdot,\cdot,\cdot)$ among trajectories with a fixed history. 
Hence, by compactness, this problem has at least one solution, see \cite[Theorem~1]{piccoli98}.

Let $(\tau^\varepsilon, q^\varepsilon(\cdot), x^\varepsilon(\cdot))$ be a solution of $(\mathrm{HP})^{'}_\varepsilon$, then it is also a solution of $\mathrm{(HP)}_\varepsilon$.
We apply Lemma~\ref{lemzeno}, and we consider the corresponding sequence $(\tau^n, q^n(\cdot), x^n(\cdot))$, which by construction has a finite number of location switchings. 
Then, by optimality and using~(\ref{convcosti}),
\begin{equation*}
\begin{split}
0 & \leq  C(\tau^\varepsilon,q^\varepsilon(\cdot), x^\varepsilon(\cdot))  - C(\tau^*,q^*(\cdot), x^*(\cdot))\\
& \leq  C(\tau^\varepsilon,q^\varepsilon(\cdot), x^\varepsilon(\cdot))- C(\tau^*,q^*(\cdot), x^*(\cdot))+\varepsilon\, \mathrm{TV}(h\circ q^\varepsilon(\cdot))\\
&\leq C(\tau^n,q^n(\cdot), x^n(\cdot))- C(\tau^*,q^*(\cdot), x^*(\cdot))+\varepsilon\, \mathrm{TV}(h\circ q^n(\cdot))\\
&\leq \mathrm{O}(\tau_\infty^*-\tau^*_{n}) + \varepsilon n |Q|,
\end{split}
\end{equation*}
where $|Q|$ is the number of locations. Choose $n = \lfloor \varepsilon^{-1/2} \rfloor$. The convergence \eqref{costizeno} follows by letting $\varepsilon$ converge to $0$.
\end{proof}

\begin{remark}\label{rk:9292}
If the rate of convergence of $\tau_n^*$ to $\tau_\infty^*$ is known, then it is possible to determine the rate of convergence in~(\ref{costizeno}). For instance, if $\tau_n^*-\tau_\infty^* \leq \mathrm{O}(n^{-\beta})$ for some $\beta >0$, then, for every $\alpha >0$, we have
$$
C(\tau^\varepsilon,q^\varepsilon(\cdot), x^\varepsilon(\cdot))- C(\tau^*,q^*(\cdot), x^*(\cdot)) \leq \mathrm{O}\left(\varepsilon^{\min(1-\alpha, \alpha\beta)}\right).
$$
\end{remark}

\appendix

\section{Further comments on condition \condition}\label{sec:remarks}
The relation between  condition \condition\ and small-time local controllability depends the continuity of the time-optimal map $\Upsilon$. 
Recall that $\Upsilon(y)$ is the minimal time needed to steer the control system \eqref{dinamica} from $y$ to $0$. 

Note that, in Definition~\ref{def:omega}, $(\Omega_{2})$ does not imply $(\Omega_{1})$ in general, as the following example shows.

\begin{example}\label{ex:marco}
Consider the control system 
$$
\dot x = u f_1(x) + v f_2(x), \quad x = (x_{1},x_{2}) \in \R^2,\quad (u,v) \in [-1,1]^{2},
$$
where $f_1(x)=\partial_{x_1}$, $f_2(x)=h(x_1) \partial_{x_2}$, with 
$$
h(x_1)
=\begin{cases}
0,& \mbox{ if } x_{1} \in [-1,1],\\
1,& \mbox{ if } x_{1} \notin (-2,2),
\end{cases}
$$
and $h$ is a smooth function with $h(x_{1}) \in [0,1]$ for every $x_{1} \in \R$.

The control system is clearly not STLC at $0$. However every point of $\R^{2}$ can be steered to the $0$ with at most two switches.
Moreover every point $y$ in the open strip $\mathcal{N} = (-1,1)\times\R$ can be steered to $0$ with two switches in time $\tau_y \leq 4 \Upsilon(y) $. Indeed consider for instance $y =(y_1,y_2)\in \mathcal{N}$ with $y_1 \geq 0, y_2>0$ (the other cases can be treated similarly). The control
$$
(u(t),v(t)) = 
\begin{cases}
(1,0), & t \in [0,2-y_1)\\
(0,-1), & t \in [2-y_1, 2-y_1+y_2)\\
(-1,0), & t \in [2-y_1+y_2, 4-y_1+y_2]
\end{cases}
$$
steers $y$ to $0$ in time $\tau_y=4-y_1+y_2 \leq 4+y_2$ while $\Upsilon(y) \geq 1+y_2$. 
Hence the control system satisfies condition $(\Omega_{2})$.
\end{example}

In the example above, the time-optimal map is not continuous at $0$. Indeed $\Upsilon(0)=0$ while $\Upsilon((0,x_{2})) \geq 2$ for every $x_{2} \neq 0$. A relationship between $(\Omega_{2})$ and $(\Omega_{1})$ can be established depending on the continuity of the time-optimal map $\Upsilon$.

\begin{proposition}
The following conditions are equivalent.
\begin{itemize}
\item[$(a)$] $\Upsilon$ is continuous at $0$.
\item[$(b)$] Condition $(\Omega_{2})$ implies condition $(\Omega_{1})$.
\end{itemize}
\end{proposition}
\begin{proof}
\noindent{$(a)\Rightarrow (b)$.} A stronger assertion actually holds, namely, $(a)$ implies $(\Omega_{1})$. Indeed, if $\Upsilon$ is continuous at $0$, since  $\Upsilon(0)=0$, for every $\varepsilon>0$, $\Upsilon^{-1}([0,\varepsilon))$ is a neighborhood of $0$ and every point in $\Upsilon^{-1}([0,\varepsilon))$ can be steered to $0$ in time less than $\varepsilon$ for the control system \eqref{dinamica}.

\noindent{$(b)\Rightarrow (a)$.} 
This is a consequence of the classical fact that, if \eqref{dinamica} is STLC at $0$, then $\Upsilon$ is continuous at $0$ (see~\cite[Theorem 2.2]{stefani-holder}).
\end{proof}

In the rest of this section, we present sufficient conditions ensuring \condition. First, note that in the simple case of a driftless control-affine system, \condition\ is a consequence of the Lie Algebra Rank Condition. In this case the number of switchings needed to reach any point in a small neighborhood of $0$ depends only on the step of the Lie algebra $\mathrm{Lie}(f_{1},\ldots,f_{m})$ at $0$. 

\begin{proposition}\label{prop:1}
For a driftless control-affine system $\dot x = \sum_{i=1}^{m} u_{i}f_{i}(x)$, if $\mathrm{Lie}_{0}(f_{1},\ldots,f_{m}) = \R^{N}$, then \condition\ is satisfied at $0$.
\end{proposition}

For control-affine systems with a drift, a sufficient condition comes from the classical result by Sussmann~\cite{suss-bound} in the single-input case. The main assumption in~\cite{suss-bound} (denoted by $(\Delta)$ in this reference) involves Lie brackets between the drift vector field and the controlled vector field (we also refer to~\cite{sharon} for more precise estimates on the number of switchings in a particular case). More precisely we have the following result.

\begin{proposition}\label{prop:2}
Consider the single-input control-affine system $\dot x = f(x)+ ug(x)$, where $f$ and $g$ are analytic vector fields in $\R^N$. If the condition $(\Delta)$ of \cite{suss-bound} is satisfied, and if the control system is STLC at $0$, then \condition\ holds true at $0$.
\end{proposition}

\begin{proof}
By~\cite{suss-bound}, the system satisfies the bang-bang property with bounds on the number of switchings (BBBNS). More precisely, for every $K$ compact and for every $T>0$, there exists $n_0\in\N^*$ such that, if $x(\cdot)$ is a time-optimal trajectory that is entirely contained in $K$ and steers the control system from $x \in K$ to $y\in K$, then there exists a time-optimal trajectory steering as well the control system from $x$ to $y$, which is moreover bang-bang with at most $n_0$ switchings, with $n_0$ depending on $K$ and $T$. 
Since the control system is STLC at $0$, the set $K= \{x\mid \Upsilon(x) \leq 1\}$ is a compact set containing $0$ in its interior. Every $x\in K$ can be steered to $0$ in time $\Upsilon(x)$ with at most $n_0$ switchings. 
\end{proof}

Linear autonomous systems generically satisfy \condition, as established next.

\begin{proposition}\label{prop:3}
If the linear autonomous control system $\dot x = A x+ B u$ satisfies the Kalman condition, then \condition\ holds true.
\end{proposition}

\begin{proof}
It suffices to write the system in Brunowsky form (see, e.g., \cite[Theorem~8, Section 4.2]{sontag}). The time-optimal control of a cascade system has a number of switchings depending only on \emph{Kronecker indices} (or \emph{controllability indices}) of the system (see also \cite{LM}).
\end{proof}

As a consequence, we have the following sufficient condition for control-affine systems.

\begin{proposition}
Consider the control affine system $\dot x = f(x) + \sum_{i=1}^{m} u_{i}g_{i}(x)$. We set
$$
G_{i} = \spann\{\mathrm{ad}^{k}_{f}g_{j}\ \mid\ 0\leq k \leq i,\ 1\leq j\leq m\}.
$$
Assume that:
\begin{itemize}
\item[$(i)$] for every $1\leq i\leq N-1$, the distribution $G_{i}$ has constant dimension near $0$;
\item[$(ii)$] the distribution $G_{N-1}$ has dimension $N$;
\item[$(iii)$] for every $1\leq i\leq N-2$, the distribution $G_{i}$ is involutive.
\end{itemize}
Then \condition\ holds true at $0$.
\end{proposition}

\begin{proof}
The result follows from Proposition~\ref{prop:3} and from the fact that the \emph{State Space Exact Linearization Problem} is solvable (see, e.g., \cite[Theorem 5.2.3]{isidori}).
\end{proof}

\section{An existence result}\label{sec:app}
For every $\varepsilon\geq 0$, consider the optimal control problem
\begin{equation*} 
\begin{cases}
&\displaystyle{\min_{u\in{\cal U}}  \left(\int_0^{t(u)}L(s,x(s),u(s))\,ds + \varepsilon\, \mathrm{TV}(u)\right),} \\
&\dot x=f(t,x,u),\quad u\in{\cal U},\\
& x(t)\in\con, \quad t\in[0,t(u)],\\ 
&x(0)\in M_{0},\quad x(t(u))\in M_{1},
\end{cases}
\eqno{(\mathrm{OCPS})_{{\varepsilon}}}
\end{equation*}
where
\begin{itemize}
\item $f:\R\times\R^N\times{\bf U}\to\R^N$ is measurable w.r.t. $t$, locally Lipschitz w.r.t. $x$,
\item $L\in\con^{0}(\R\times\R^{N}\times\R^{m})$, 
\item $\con=\{x\in \R^{N}\mid h_1(x)\geq 0,\dots, h_l(x)\geq 0\}$ for some $h_1,\dots h_l\in\con^0(\R^N)$,
\item ${\bf U}\subset\R^m$ is compact, 
\item $M_{0}$ and $M_{1}$ are compact subsets of $\con$.
\end{itemize}
Here, $\cal U$ is still defined by \eqref{defcalU}.

\begin{theorem}\label{exth}
Assume that:
\begin{itemize}
\item[(i)] there exists $\bar u\in{\cal U}$ having bounded variation, steering the control system $\dot x=f(t,x,u)$ from $M_0$ to $M_1$, and whose corresponding trajectory satisfies the state constraint $x(t)\in\con$, for every $t\in[0,t(\bar u)]$;
\item[(ii)] there exists $b>0$ such that, for every $u\in{\cal U}$ steering the control system from $M_0$ to $M_1$, its corresponding trajectory $x_u$ satisfies $t(u)+\Vert x_u(\cdot)\Vert_{\infty}\leq b$.
\end{itemize}
Then, for every $\varepsilon>0$, the optimal control problem $(\mathrm{OCPS})_{{\varepsilon}}$ has at least one solution.
\end{theorem}

Note that existence is not ensured for $\varepsilon=0$. The fact that $\varepsilon>0$ is crucial here. The difference with usual existence theorems is that, in the proof below, we use in an instrumental way the total variation term. Note the remarkable fact that, in contrast to usual existence theorems (see \cite{cesari-book}), we do not assume, here, that the set of extended velocities \eqref{defV} is convex. This classical assumption can be removed thanks to the use of the total variation term.

\begin{proof}
The proof follows the lines of~\cite[Theorem 5.14 and 6.15]{trelat_book}, with an adaptation to the bounded variation context.
Let
\bqnn
\delta=\inf\left(\int_0^{t(u)}L(s,x(s),u(s))\,ds + \varepsilon\, \mathrm{TV}(u)\right),
\eqnn
where the infimum is taken among all controls $u\in{\cal U}$ steering the control system from $M_0$ to $M_1$ and whose corresponding trajectory satisfies the state constraint $x(t)\in\con$, for every $t\in[0,t(u)]$. 
Let $x_n(\cdot)$ be a sequence of admissible trajectories, corresponding to a minimizing sequence of admissible controls $u_n:[0,t(u_n)]\to {\bf U}$, i.e., 
\bqnn
\lim_{n\to\infty}\left(\int_0^{t(u_n)}L(s,x_n(s),u_n(s))\,ds + \varepsilon\, \mathrm{TV}(u_n)\right)=\delta .
\eqnn
Using Assumptions~$(i)$ and $(ii)$, for $n$ sufficiently large we have
\bqnn
\varepsilon\, \mathrm{TV}(u_n)\leq \int_0^{t_(\bar u)}L(s,x_{\bar u}(s),\bar u(s))\,ds + \varepsilon \,\mathrm{TV}(\bar u) + C,
\eqnn
for some constant $C\geq 0$, and since $t(u_n)$ is bounded by $b$, extending $u_n$ by $0$ for $t>t(u_n)$, we infer that the sequence $(u_n)_{n\in\N}$ is bounded in the set $\mathrm{BV}([0,b],\R^m)$ of bounded variation functions from $[0,b]$ to $\R^m$. Since the embedding $\mathrm{BV}([0,b],\R^m)\hookrightarrow L^1([0,b],\R^m)$ is compact (see \cite{EvansGariepy}), up to some subsequence, $(u_n)_{n\in\N}$, converges to some $u_{\varepsilon}\in L^1([0,b],\R^m)$ for the strong topology of $L^1$. Still up to some subsequence, $x_n(0)$ converge to some $x_{\varepsilon}^0\in\R^N$, $u_n$ converges to $u_{\varepsilon}$ almost everywhere and $t(u_n)$ converges to $t(u_\varepsilon)$, and thus $u_{\varepsilon}:[0,t(u_\varepsilon)]\rightarrow {\bf U}$ takes values in ${\bf U}$.

Let us prove that $u_{\varepsilon}:[0,t(u_\varepsilon)]\to {\bf U}$ is a solution of $(\mathrm{OCPS})_{{\varepsilon}}$. By a standard Gronwall argument (see \cite[Theorem~1 p.~56]{sontag}, or see \cite{Trelat2000,trelat_book}), the convergence almost everywhere of $u_n$ to $u_{\varepsilon}$ implies that $x_n(\cdot)$ converges uniformly to $x_{\varepsilon}(\cdot)$, where  $x_{\varepsilon}(\cdot)$ is the trajectory corresponding to the control $u_{{\varepsilon}}$ and starting at $x^0_{{\varepsilon}}$.
In particular, we get that $x_{\varepsilon}(t)\in\con$ for every $t\in[0,t(u_\varepsilon)]$ and, by compactness of $M_0$ and $M_1$, we obtain that $x_{\varepsilon}(t(u_\varepsilon))\in M_1$. Hence $u_{\varepsilon}$ is an admissible control. Moreover, $L(t,x_n(t),u_n(t))$ converges to $L(t, x_{\varepsilon}(t), u_{\varepsilon}(t))$ for almost every $t$.
Hence, using Assumption $(ii)$ and the dominated convergence theorem, we conclude that 
\begin{align}
\lim_{n\to\infty}\int_0^{t(u_n)} & L(t,x_n(t),u_n(t))\,dt\nonumber\\
& =\int_0^{t(u_\varepsilon)}L(t, x_{\varepsilon}(t), u_{\varepsilon}(t))\,dt . \label{cint}
\end{align}
On the other hand, by lower semicontinuity of the functional $\mathrm{TV}(\cdot)$, we have
\begin{equation}\label{ctv}
\mathrm{TV}(u_{\varepsilon})\leq \liminf_{n\to\infty} \mathrm{TV}(u_n).
\end{equation}
Using \eqref{cint}, \eqref{ctv} and since $u_n$ is a minimizing sequence, we infer that
$\int_0^{t(u_\varepsilon)}L(t, x_{\varepsilon}(t), u_{\varepsilon}(t))\,dt + {\varepsilon} \mathrm{TV}(u_{\varepsilon})\leq\delta$, which implies that $u_{\varepsilon}$ is optimal.
\end{proof}

\section*{Acknowledgement}
The authors are grateful to Jean-Michel Coron, Matthias Kawski, Andrey Sarychev, and Mikhail Il'ich Zelikin for useful discussions and valuable suggestions. 
The last author acknowledges the support by FA9550-14-1-0214 of the EOARD-AFOSR.

\bibliographystyle{IEEEtran}
\bibliography{IEEEabrv,biblio_fuller}

\end{document}